# Too Many Hats

Rob Pratt, Stan Wagon, Michael Wiener, and Piotr Zielinski

**Abstract.** A puzzle about prisoners trying to identify the color of a hat on their head leads to a version where there are $k$ more hats than prisoners. This generalized puzzle is related to the independence number of the arrangement graph $A_{m,n}$ and to Steiner systems and other designs. A natural conjecture is that perfect hat-guessing strategies exist in all cases, where "perfect" means that the success probability is $1/(k + 1)$. This is true when $k = 1$, but we show that it is false when $k = 2$. Further, we present a strategy with success rate at least $1/O(k \log k)$, independent of the number of prisoners.

Recent years have seen a profusion of hat puzzles, which seek strategies that some prisoners can use to gain their freedom. These riddles are attractive because there are often strategies that have much higher success rates than one would think possible. Here is the classic case:

> There are $n$ prisoners, Alice, Bob, Charlie,… under the care of a warden who lines them up in order. He has $n$ hats—each is either red or blue—and randomly places one hat on each prisoner's head. Each prisoner can see the hats of only the prisoners in front of them: Alice sees all but her own, Bob sees $n - 2$ hats, and so on; the last prisoner sees nothing. Alice then guesses her color; all prisoners can hear the guesses. Then Bob guesses his hat color, and so on for all prisoners. If they all are correct, they will all be freed. Otherwise, none are freed.
>
> The prisoners know the rules and can devise a strategy in advance; no communication other than the guesses is allowed once the hats are placed. Find their best strategy.

Simply guessing randomly gives the prisoners a $2^{-n}$ chance of getting them all correct. Alice is in a bad position, because for her the probability of success cannot exceed 50%: she is guessing from two equally likely choices. But the chance of correctness for everyone else can be increased to absolute certainty if they use a clever parity strategy. Alice will declare "red" if the number of reds she sees is even; otherwise she says "blue". Hearing this, Bob can easily determine, from the colors he sees, whether his hat is red or blue, and the same for all the other prisoners.

There are many variations on this folklore puzzle, several of which are discussed in [6] (see also [3] for a version where each prisoner sees all other prisoners). For example, there might be three hat colors for the $n$ prisoners. Or there might even be infinitely many prisoners, a topic discussed in detail in the book by Hardin and Taylor [5]. Here we study the variation where the hats all have distinct colors. The case of $n$ prisoners and $n$ differently colored hats is uninteresting, as Alice can immediately deduce her color. But Tanya Khovanova [6] presented the case where there are $n + 1$ hat colors and the warden just discards the unused hat, a nice variation due to K. Knop and A. Shapovalov.

For the extra-hat case, Alice must make a choice from the two hats she does not see, and so she can succeed at most half the time. But again the other prisoners can be certain of their color if they all agree to a parity strategy. The twist is that it is the parity of a permutation that they must analyze. The prisoners, numbered 1 (Alice) through $n$, imagine a ghost prisoner, number $n + 1$, who wears the missing hat. The colors are identified with $1, 2, …, n + 1$ and the prisoners assume that the permutation $\pi$ of $1\,2\,…\,(n + 1)$ induced by the hats is an even permutation. That assumption gives each prisoner only one possibility for his or her declaration. So they succeed when $\pi$ is in fact even; half the permutations are even, so they succeed half the time. This is a perfect situation, because the strategy wins 50% of the time and this is best possible.

## Two Extra Hats

The preceding puzzle leads naturally to the case of two or more extra hats, and there are some surprises as well as interesting connections to graph theory, Steiner systems, Latin squares, and ordered designs. We use $n$ for the number of prisoners and $k$ for the number of extra hats. Because Alice must choose from $k + 1$ hats (she sees $n - 1$ of the $n + k$ hats), the chance of any strategy's success is never greater than $1/(k + 1)$. We consider here only deterministic strategies (as opposed to probabilistic ones). Call a strategy *perfect* if its success probability is $1/(k + 1)$. A natural conjecture is that a perfect strategy exists in all cases.

We start with $k = 2$ and introduce two ghosts (numbered $n + 1, n + 2$) who wear the unused hats. If there are only two prisoners, it is easy to find a perfect strategy. The prisoners will assume that the hat assignment is $(1, 2)$, $(2, 3)$, $(3, 4)$, or $(4, 1)$. Equivalently, Alice subtracts 1 from the color she sees, while Bob adds 1 to what he hears, working modulo 4. Because there are 12 possible hat states (in general, $(n + k)!/k!$ states), the success rate is $4/12$ or $1/3$. For more than two prisoners, the problem becomes complicated and there are several types of strategies. We start with a natural arithmetic strategy, leaving the details as an exercise.

**Modular Arithmetic Strategy (Larry Carter).** The prisoners assume that the sum of the hats on their heads is 1 (mod $n + 2$) (if $n \equiv 2 \pmod 4$) and 0 (mod $n + 2$) otherwise. We omit the details, but the probability that the modular assumption is correct is $1/(2 \lceil n/2 \rceil + 1)$ and, when it is correct, they win. For two prisoners this strategy has success probability $1/3$; for three prisoners, it is $1/5$.

Because the hat-sums are not equidistributed modulo $n + 2$, using a single residue such as 0 is not optimal. While the use of the two residues 0 and 1 maximizes the success probability for this strategy, the modular arithmetic strategy is far from optimal when $n > 2$.

The next strategy uses parity for both numbers and permutations; it shows that the prisoners, regardless of how many there are, always have at least a $1/4$ chance of winning. When dealing with the parity of the full permutation of the colors, we must make some assumption about the order of the ghost colors. Because the ghosts can exchange hats at will, any assumption is allowed. For the following strategy we assume that an even ghost color always precedes an odd. This turns any hat assignment into $\pi$, a permutation of $\{1, 2,\ldots, n, n + 1, n + 2\}$.

**Double Parity Strategy.** The prisoners assume that the hat assignment satisfies:

1. The unused colors have different parity.

2. The permutation $\pi$ is even; here the definition of $\pi$ assumes that the unused colors are in the order (even, odd).

If the unused colors all have the same parity and Alice, or another prisoner, can deduce that, then this strategy is undefined, and the prisoners lose.

*The strategy wins when the assumption holds.* Consider Alice, who, because of what she sees, knows the three missing colors. By (1), her color must have the parity that appears twice among the three and so there are two choices. Suppose they are both even; then they appear in positions 1 and $n + 1$, and there is only one possibility that leads to $\pi$ being even. The odd case is similar. Argue the same way for Bob, because he knows that Alice's declaration is correct, and so on inductively for all the prisoners. □

*The success probability is at least* $1/4$. Suppose $n$ is even. There are $(n + 2)/2$ choices for the odd unused color, and same for the even. The used colors can be permuted in any way so that the final permutation is even; this is half of the permutations. So the winning count is $\frac{1}{2} n! \left(\frac{n}{2} + 1\right)^2$. The corresponding probability is the ratio of this to $(n + 2)!/2$, or $\frac{1}{4} + \frac{1}{4(n+1)}$. The odd case is similar, yielding $\frac{1}{4} + \frac{1}{4(n+2)}$. The limiting probability is $1/4$. □

This asymptotic success rate of 25% is the best such result we know of, but far from the $1/3$ that perfect strategies attain. For three prisoners the hat assignments in the assumed set of the double parity strategy are

```
123  134  145  152  215  231  253  312  325
341  354  413  435  451  514  521  532  543
```

with size 18 and success probability $18/60$, or $3/10$. We shall see in a moment that a perfect strategy exists when $n = 3$. But for large $n$, the double parity strategy is the best known.

In many cases there are perfect strategies, which outperform the preceding ones. A perfect strategy for $n$ prisoners and $k$ extra hats has an important connection to the arrangement graph $A_{n+k,n}$: this graph has as vertices all ordered $n$-tuples consisting of distinct integers chosen from 1 through $n + k$, with two vertices being adjacent if the corresponding tuples differ in exactly one position. Thus the vertex set consists of all possible hat assignments. We use $\alpha_{n,k}$ to denote the independence number (size of largest independent set) of $A_{n+k,n}$.

**Theorem 1.** For $n$ prisoners and $k$ extra hats, a perfect strategy exists if and only if $\alpha_{n,k} = (n + k)!/(k + 1)!$ (i.e., there is an independent set in $A_{n+k,n}$ having size that is $1/(k + 1)$ of the vertex count).

*Proof.* If $X$ is an independent set, the prisoners can assume that the hat assignment lies in $X$; if it does in fact do so, the color of each prisoner's hat is uniquely determined and the prisoners will win; if $X$'s size is $1/(k + 1)$ of the vertex count then the resulting strategy is perfect. Conversely, any strategy leads to the set of all hat assignments for which the strategy wins; this set is an independent set in $A_{n+k,n}$ because an edge in this set would mean that one prisoner's color is not uniquely determined. If the strategy is perfect, then the size of the independent set is as claimed. □

An independent set as in Theorem 1 is called a *perfect independent set*. Figure 1 shows $A_{6,2}$ (the 30 vertices are ordered pairs from 1 through 6; edges are *all* vertical and horizontal connections, not just the nearby ones shown; e.g., $(5, 1) \leftrightarrow (5, 4)$). A perfect independent set of size six is shown; if the hat assignment is one of these six, then, if Alice sees color $i$ she knows her color is $i - 1$ (mod 6), and the same for Bob with $i + 1$. The larger the independent set, the higher the probability of success and so the best possible strategy requires computing $\alpha_{n+k,n}$; we use $V_{n+k,n}$ for the vertex count of $A_{n+k,n}$, which is $(n + k)!/k!$. The graph $A_{n+1,n}$ (the case $k = 1$) is bipartite, with parts defined by the parity of the permutation $\pi$. This gives $\alpha_{n+1,n} = V_{n+1,n}/2$ and so yields a perfect strategy, identical to the one-hat-too-many solution given earlier. When $k = 2$ this graph is a Cayley graph of the alternating group graph $AG_{n+2}$ [11]; this graph has vertices for every group element and an edge connecting permutations that differ by either $(1, i, 2)$ or $(1, 2, i)$, where these are in cycle notation, $i \geq 3$, and the first is used if $i$ is odd and the second if $i$ is even.

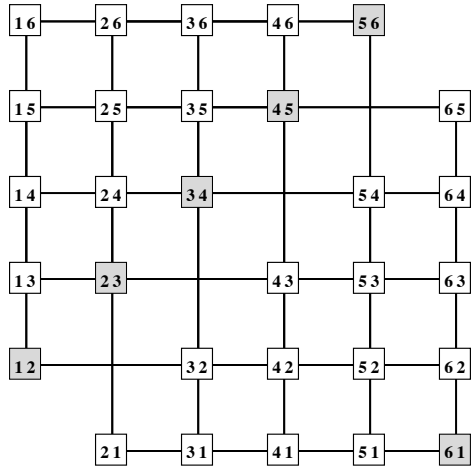

**Figure 1.** The arrangement graph $A_{n+k,n}$, where $n = 2$ and $k = 4$ and with vertex coordinates equal to the 2-tuple. The edges are all horizontal or vertical connections. The gray vertices are a perfect independent set and give a perfect strategy.

Further, a perfect independent set in $A_{n+k,n}$ is exactly an *ordered design* $OD_1(n - 1, n, n + k)$ (see [1]). An $OD_1(t, n, v)$ is an $n \times \binom{v}{t} t!$ array with entries from 1 to $v$ so that

- each column has $n$ distinct entries;
- each collection of $t$ rows contains each possible $t$-tuple exactly once among its columns.

For example, the following $3 \times 20$ array is an $OD_1(2, 3, 5)$. There are $2\binom{5}{2} = 20$ ordered pairs from $\{1, 2, 3, 4, 5\}$

and for each way of selecting two rows, all 20 pairs appear as columns.

```
1 1 1 1 2 2 2 2 3 3 3 3 4 4 4 4 5 5 5 5
2 3 4 5 1 3 4 5 1 2 4 5 1 2 3 5 1 2 3 4
4 2 5 3 4 5 3 1 2 5 1 4 5 3 1 2 3 1 4 2
```

Th existence of an $OD_1(n - 1, n, n + k)$ is equivalent to the existence of a perfect independent set in $A_{n+k,n}$. Each column of the array corresponds to a vertex in $A_{n+k,n}$; the number of columns is $\binom{n+k}{n-1}(n-1)!$, which is the same as $\frac{(n+k)!}{(k+1)!}$, the size of a perfect independent set, and two columns cannot differ in one position only, as that would violate the last point in the definition. Thus the existence of such designs lead to perfect hat strategies (and the converse is also true). The cases $n \leq 4$ of the next result were known to researchers in ordered designs, but the results for $n \geq 5$ are new.

**Theorem 2.** When $k = 2$, perfect strategies exist for $n \leq 6$ and for no other values.

*Proof.* $n = 2$: $\alpha_{4,2} = 4$ by $(1, 2), (2, 3), (3, 4), (4, 1)$, a perfect independent set because $12/3 = 4$.

$n = 3$: If the prisoners assume that the hat-color vector is in $S = \{(a, b, (3a + 3b)_{\text{mod } 5}) : 1 \leq a, b \leq 5, \ b \neq a\}$, then they win whenever the assumption is correct because any two of $a, b, 3(a + b)$ determine the third, and because $3a + 3b \not\equiv a$ or $b$ (mod 5). Note that $S$ is invariant under permutation of the first two elements. This set $S$ is identical to the columns of the earlier example of an $OD_1(2, 3, 5)$. The 20 triples—hat assignments—show that $\alpha_{5,3} = 20$ and $S$ is a perfect independent set because 20 is one third of $V_{5,3} = 60$.

$n = 4$: Consider this matrix, derived from Teirlinck's work [9] on ordered designs:

$$M = \begin{pmatrix} 5 & 3 & 1 & 4 & 2 & 0 \\ 3 & 5 & 4 & 2 & 0 & 1 \\ 1 & 4 & 5 & 0 & 3 & 2 \\ 4 & 2 & 0 & 5 & 1 & 3 \\ 2 & 0 & 3 & 1 & 5 & 4 \\ 0 & 1 & 2 & 3 & 4 & 5 \end{pmatrix}.$$

Note that $M$ is a symmetric Latin square with constant diagonal. The method of construction of $M$ is given in the $n = 4$ case in the next section. Let $S$ be the set of 4-vectors $(a, b, c, d)$ from $\{1, 2, 3, 4, 5, 6\}$ so that $M_{a,b} = M_{c,d}$. We claim that $S$ is a perfect independent set. Two vectors in $S$ cannot differ in exactly one coordinate because each row (and column) of $M$ has distinct entries. The common value in the defining equation can be any of 0 through 4; assume it is 0. We claim there are 24 vectors in $S$ for this value. Suppose $(a, b, c, d) \in S$. Then, because $M$ is symmetric, $S$ contains $(a, b, d, c), (b, a, c, d), (b, a, d, c), (c, d, a, b), (d, c, a, b), (c, d, b, a)$, and $(d, c, b, a)$. That is, $S$ is invariant under the eight-element group $G$ generated by $(1\,2)$ and $(1\,4)(2\,3)$. And there are three possibilities for the orbit-generators: 1625, 1634, and 2534; these correspond to the three 0s in the upper right quadrant. So there are 24 in all, as claimed. This count works for any entry in place of 0, so $|S| = 5 \cdot 24 = 120$ and $S$ is perfect. To see the vectors in $S$ transform each of the following 15 vectors (representatives of the $G$-orbits) by the eight permutations in $G$:

1235 1246 1326 1345 1423 1456 1524 1536 1625 1634 2356 2436 2534 2645 3546

$n = 5$: Some computer searching led us to a perfect strategy of 840 vectors, useful tricks being to enforce some symmetry, or to assume that the perfect $n = 4$ strategy embeds in the one being sought. One approach is to use ILP (integer linear programming) with 0-1 variables $x_{A,i}$ where $A$ is a possible $(n - 1)$-vector and $i$ a possible last entry in the $n$-vector. When this is seeded by using the result of appending 7 to each vector in an ILP-based solution for $n = 4$, it takes only a half-second to find the desired 840-sized perfect independent set for $n = 5$. The full set can be described thus: extend each of the following 42 hat assignments by the order-20 group generated by $(2\,4)(3\,5)$ and $(1\,5\,3\,2)$.

```
12563  12643  12354  12534  13654  23564  12465  12635  13645  12346  12456  23546  41327  42317
51247  51427  52137  52347  53127  53417  54237  54317  61237  61347  61457  61527  62147  62357
62417  62537  63157  63217  63427  63547  64137  64257  64327  64517  65127  65247  65317  65437
```

$n = 6$: The same ideas used for $n = 5$ work. The set we found is invariant under the group of order 120 generated

by (1 2) (4 5) and (2 6 3 5). The group orbit of the following 56 elements gives the perfect independent set of size 6720.

```
123456  123468  123475  123487  123548  123567  123574  123586
123645  123657  123678  123684  123746  123758  123764  123785
123847  123854  123865  123876  124568  124576  124658  124687
124765  124786  124857  124875  125768  125867  134578  134587
134657  134685  134756  134768  134865  134876  135678  135786
145687  145876  234567  234586  234678  234687  234758  234765
234856  234875  235867  235876  245768  245786  345678  345867
```

When $n = 7$, the double parity strategy leads to 50400 vectors, but this can be improved to 50880 (hence a 28% chance of success) by an ILP approach that assumes the set is invariant under the 120 permutations of indices generated by (12) (37) and (2654) (37). A perfect strategy requires 60480 vectors; it came as a surprise when a backtracking search showed that such a 60480-sized set does not exist (see Table 1), thus disproving the conjecture that perfect strategies always exist when $k = 2$. This result means that a perfect strategy does not exist for $n \geq 7$, because such a strategy for $n$ easily leads to one for $n - 1$ (delete $n + 2$ from all vectors ending in $n + 2$). □

|  | $n = 1$ | $n = 2$ | $n = 3$ | $n = 4$ | $n = 5$ | $n = 6$ | $n = 7$ |
| --- | --- | --- | --- | --- | --- | --- | --- |
| $k = 1$; bipartite | 1 | 3 | 12 | 60 | 360 | 2520 | 20 160 |
| $k = 2$ | 1 | 4 | 20 | 120 | **840** | **6720** | [50880, 60479] |
| $k = 3$ | 1 | 5 | 30 | [204, 206] | <1648 | <14832 | <148320 |
| $k = 4$ | 1 | 6 | 42 | 336 | **3024** |  |  |
| $k = 5$ | 1 | 7 | 56 | 504 |  |  |  |
| $k = 6$ | 1 | 8 | 72 | 720 | 7920 | 95 040 |  |
| $k$ | all perfect | all perfect | all perfect | all perfect except $k = 3$ | ∞ many perfect | ∞ many perfect | ∞ many perfect |

**Table 1.** The independence number $\alpha_{n+k,n}$. Black entries are $(n + k)!/(k + 1)!$ and indicate perfect strategies; the black non-bold entries follow from known results about ordered designs; bold entries are new. Red intervals indicate new results that give bounds on $\alpha_{n+k,n}$.

## More Hats

The problem can be studied when there are three or more unused hats. The double parity strategy for $k = 2$ extends to show that for $k \geq 2$ there is a strategy for $n$ prisoners having success probability greater than $1/(e\,k^2)$, independent of $n$. For this extension, we again imagine that ghosts wear the unused hats, and that the ghosts are in the order specified in condition (2) below. Let $t = \lfloor k^2/2 \rfloor$. Then the prisoners make the following three assumptions.

1. The unused colors are distinct modulo $t$.

2. The permutation $\pi$ is even, where the unused colors are assumed to be in the order of their mod-$t$ residues.

3. The mod-$t$ sum of the unused colors is $\sigma$ (where $\sigma$ is chosen to maximize the success rate).

Note that (3) follows from (1) when $k = 2$ and $t = 2$ ($\sigma$ being 1). The proof that this assumption's truth leads to a win is the same as for the $k = 2$ case discussed earlier: (1) and (3) narrow Alice's possibilities to one or two colors; if two, (1) and (2) yield the correct color. Calculating the probability that the assumption holds for the hat assignment requires a little work. The key is to first study the probability of (1), a problem identical to the classic birthday puzzle (with $t$ days and $k$ people). We omit the details, but the probability that (1) holds is at least $t!/\bigl((t - k)!\,t^k\bigr)$; standard factorial approximations and bounds show that this is at least $1/e$ for our choice of $t$. Now an averaging argument implies that there is some $\sigma$ so that the probability of (1) and (2) is at least $1/(e\,t)$, yielding the lower bound $2/(e\,k^2)$. Condition (2) then reduces this to $1/(e\,k^2)$. The details of this analysis show why, in the choice of $t$, the exponent 2 and coefficient $1/2$ are the best choice, in the limit as $k \to \infty$. But when $2 \leq k \leq 5$, a better choice is $t = k$ with the target sum $\sigma = 1, 0, 2, 0$, respectively for the four mod-4 cases. The asymptotic probabilities are then $k!/(2\,k^k)$, which is better than $1/(e\,k^2)$.

Having a reciprocal quadratic success rate is nice (and it can be improved: see next section) but that is far from a perfect strategy. We can get perfect strategies (same as perfect independent sets in $A_{n+k,n}$) in many cases by focusing on $n$.

When $n = 2$ the cyclic method used earlier gives a perfect independent set for all $k$; $\alpha_{k+2,2} = k + 2$ via $\{(a, (a + 1)_{\mod k+2}) : 1 \le a \le k + 2\}$.

When $n = 3$ and $k$ is even, then $k + 3$ is odd and the set $\{a, b, ((a + b)(k + 4)/2)_{\mod k+3}\}$ (where $1 \le a, b \le k + 3$; $a \ne b$) is a perfect independent set; the verification uses the fact that $(k + 4)/2$ is the mod-$(k + 3)$ inverse of 2. When $k$ is odd use $\{a, b, M_{a,b}\}$ where $1 \le a, b \le k + 3$, $a \ne b$, and $M$ is the $(k + 3) \times (k + 3)$ matrix forming a Latin square with entries from 1 through $k + 3$ that is *idempotent* (meaning: $M_{i,i} = i$ for all $i$). A simple construction for $M$ is given in [2, p. 36]. Here is a description (we omit the proof that it is an idempotent Latin square). Let $m = k + 3$ and $h = m/2$. Start with $1, 2,\ldots, m$ down the main diagonal and then place $m$ in the diagonal above the main one, except use the leftmost slot in the penultimate row. Then fill the bottom and right borders by rotated versions of $1,\ldots, m - 1$ with 1 at position $\frac{h}{2} + 1$ in the bottom row and at position $h/2$ in the rightmost column. Fill the remaining space by extending each diagonal entry along the corresponding back-diagonals; more precisely, each such remaining square gets $h(i + j) \pmod m$. When $n = k = 3$, $M$ is this $6 \times 6$ matrix:

$$\begin{pmatrix} 1 & 6 & 2 & 5 & 3 & 4 \\ 4 & 2 & 6 & 3 & 1 & 5 \\ 2 & 5 & 3 & 6 & 4 & 1 \\ 5 & 3 & 1 & 4 & 6 & 2 \\ 6 & 1 & 4 & 2 & 5 & 3 \\ 3 & 4 & 5 & 1 & 2 & 6 \end{pmatrix}.$$

When $n = 4$ perfect strategies exist for any $k$, except $k = 3$, as proved by Teirlinck [9, p. 370–372] (he used the language of orthogonal arrays and quasigroups). The negative result when $k = 3$ is that $\alpha_{7,4} \le 209$; this was proved by C. Colbourn. Using ILP we found an independent set of size 204, and then more computer searching eliminated 207; therefore $204 \le \alpha_{7,4} \le 206$. When $k$ is even, Teirlinck's methods yield a perfect strategy as follows, extending the method presented earlier when $k = 2$. Define the symmetric Latin square $M$ thus, where the fourth case is reduced modulo $k + 3$, with residue from $\{1, \ldots, k + 3\}$

$$M_{i,j} = \begin{cases} k + 3, & \text{if } i = j \\ j, & \text{if } i = k + 3 \\ i, & \text{if } j = k + 3 \\ \frac{1}{2}(k + 4)(i + j), & \text{otherwise.} \end{cases}$$

The last case uses uses the fact that $(k + 4)/2$ is the inverse of 2 (mod $k + 3$). Then the same proof as when $k = 2$ works; the size of $S$ is $8(k + 3)\binom{(k + 4)/2}{2}$, which simplifies to the perfect count $(k + 4)(k + 3)(k + 2)$. The case of $k$ odd is quite a bit more complicated.

When $n = 5$, the negative result for $n = 4, k = 3$ gives the same for $n = 5, k = 3$ by the method mentioned in the $k = 2, n = 7$ case. For $k = 4$, there is a perfect independent set. We used ILP to find such a set with an interesting symmetry property. Consider this set $S$ of 126 hat assignments:

```
12349  12356  12367  12378  12384  12395  12458  12463  12475  12487  12496  12564  12573  12589
12597  12679  12685  12698  12786  12794  12893  13457  13465  13472  13486  13498  13568  13579
13582  13594  13674  13689  13692  13785  13796  13897  14569  14576  14583  14592  14678  14682
14697  14789  14793  14895  15672  15687  15693  15784  15798  15896  16783  16795  16894  17892
23451  23468  23476  23485  23497  23569  23574  23587  23598  23675  23681  23694  23789  23791
23896  24567  24579  24586  24593  24671  24689  24695  24783  24798  24891  25678  25683  25691
25781  25796  25894  26784  26793  26897  27895  34562  34578  34589  34596  34679  34687  34691
34781  34795  34892  35671  35684  35697  35786  35792  35891  36782  36798  36895  37894  45673
45681  45698  45782  45791  45897  46785  46792  46893  47896  56789  56794  56892  57893  67891
```

Let $X$ be the result of permuting the first four entries in each entry of $S$ all 4! possible ways. Then $|X| = 24 \cdot 126 = 3024$ and $X$ is a perfect independent set in $A_{9,5}$ (and hence a new ordered design $OD_1(4, 5, 9)$). The existence of a perfect strategy in the case of $n = k = 5$ is an interesting open question. The result we discuss next

shows that for $n = 5$, perfect strategies exist when $k$ is one of 6, 10, 12, 16, 18, 22, 24, 26, 28, 30, 31, 36, and infinitely many others.

The following result of Teirlinck [7, p. 36] gives, for any value of $n$, infinitely many values of $k$ admitting perfect strategies. His theorem, translated from the language of ordered designs, is that a perfect strategy exists for $k$ and $n$ whenever the prime factorization $\Pi p_i^{a_i}$ of $k + 1$ satisfies $\Pi a_i(p_i - 1) \geq n$. In particular, this holds whenever $k \geq n$ and $k + 1$ is prime.

Another technique for getting perfect strategies involves Steiner systems $S(n - 1, n, m)$ (see [4, 12]). Such a system is a set of $n$-subsets of $\{1, 2, \ldots, m\}$ such that every $n - 1$ set appears exactly once in one of the $n$-sets. If a Steiner system $S(n - 1, n, n + k)$ exists then one can permute all its elements in all possible ways to get a perfect strategy. For example, $S(4, 5, m)$ exists when $m = 11$, giving a perfect strategy for $n = 5$ and $k = 6$. However, Teirlinck [8] proved that whenever a Steiner system $S(n - 1, n, n + k)$ exists, then his prime-factorization theorem just given applies to the parameters. Therefore a Steiner system cannot give a new perfect strategy. It is worth noting that the strategies from Steiner systems are stronger than the others in the sense that a prisoner need not see which of the other prisoners has which hats; he or she need see only the set. More precisely, if the rules were changed so that the prisoners see only the hat colors and cannot identify other prisoners by sight or by their voices, strategies based on Steiner systems still work.

## Conclusion

It is remarkable that a simple hat puzzle has connections to several different areas. Several intriguing open questions remain. The main question arises from the natural, but false, conjecture that the best strategy wins with reciprocal probability $k + 1$, independent of $n$. The double parity strategy achieves $1/(e\,k^2)$. We have found a strategy, based on an error-correcting code in [10], that succeeds with probability $1/O(k \log k)$ (see Appendix). But can this be improved?

*Question 1.* Is there a strategy for each $k, n$ so that the overall success rate in all cases is $1/O(k)$?

*Question 2.* Can anything more be said about the cases $n, k$ for which a perfect strategy exists? In particular, is there a perfect strategy when $n = k = 5$?

*Question 3.* Can $1/4$ can be improved as an asymptotic success probability when there are two extra hats?

Acknowledgement. We are grateful to Luc Teirlinck for helpful comments regarding ordered designs.

Rob Pratt (SAS Institute Inc., Cary, NC 27513; rob.pratt@sas.com)

Stan Wagon (Macalester College, St. Paul, MN 55105; wagon@macalester.edu)

Michael Wiener (Ottawa, Canada; michael.james.wiener@gmail.com)

Piotr Zieliński (Boston, MA; piotr.zielinski@gmail.com)

# Appendix. A strategy with success rate $1/O(k \log k)$.

Some ideas from the theory of error-correcting codes lead to a strategy that has a success rate of $1/O(k \log k)$. Throughout this discussion logarithms are to the base 2.

## The strategy

Let $x_1, \ldots, x_n$ be the prisoners' hat-colors and $y_1, \ldots, y_k$ the unused colors, in increasing order. For any $u, v \in \mathbb{N}$, let $\text{bit}(u, v)$ be the 0-based index of the highest bit on which the binary representations of $u$ and $v$ differ; for example, $\text{bit}(3, 6) = \text{bit}(0011_2, 0110_2) = 2$. Define the sequence $b_1, \ldots, b_{k-1}$ by $b_i = \text{bit}(y_i, y_{i+1})$, and set $b_{\max} = \text{bit}(1, n + k)$. Pick a positive integer $B$. The strategy $\sigma_{n,k,B}$ has the prisoners assume that the augmented hat configuration (also denoted by $\sigma$) satisfies the following.

(i) $x_1, \ldots, x_n, y_1, \ldots, y_k$ is an even permutation;

(ii) for all $i$, $b_i > b_{\max} - B$;

(iii) $\sum_i b_i \equiv 0 \pmod{B}$;

(iv) $\sum_{b_i < b_{i+1}} i \equiv 0 \pmod{k+1}$.

For $B = \lceil 2 \log k \rceil$, properties (i) and (ii) reduce the size of $\sigma_{n,k,B}$ from what it would be by (iii) and (iv) by a constant factor. Condition (i) divides the success probability by 2. The condition that $b_i > b_{\max} - B$ means that no two $y_i$ share the same top $B$ bits. For a given pair, the probability of this not being the case is $2^{-B}$, and there are $k^2/2$ such pairs, so the probability of this failing for at least one pair $(y_i, y_j)$ is at most $\frac{1}{2} 2^{-B} k^2 \leq \frac{k^{-2} k^2}{2} = \frac{1}{2}$ (assuming $B = \lceil 2 \log k \rceil$). So this condition also divides the probability of success by a constant that is at most 2.

Properties (iii) and (iv) define an error-correcting code capable of correcting a single deletion [10].

### Example: $n = 4, k = 2$.

Here $B = 2$ and the independent set $\sigma_{4,2,2}$ by the preceding algorithm is shown below; it has 96 elements. A perfect set has size 120, so this is not optimal. But the point of this construction is the asymptotic behavior.

```
125346  132546  153246  213546  235146  251346  315246  321546  352146  512346  523146  531246
123645  136245  162345  216345  231645  263145  312645  326145  361245  613245  621345  632145
124536  145236  152436  215436  241536  254136  412536  425136  451236  514236  521436  542136
126435  142635  164235  214635  246135  261435  416235  421635  462135  612435  624135  641235
135426  143526  154326  314526  345126  351426  415326  431526  453126  513426  534126  541326
134625  146325  163425  316425  341625  364125  413625  436125  461325  614325  631425  643125
234516  245316  253416  325416  342516  354216  423516  435216  452316  524316  532416  543216
236415  243615  264315  324615  346215  362415  426315  432615  463215  623415  634215  642315
```

The choice of the residue 0 in (iii) and (iv) might not be the best and one can use any target residues in place of 0. Then an averaging argument shows that there is at least one choice of the two target residues so that the probability that (iii) and (iv) hold is at least $1/(B(k+1))$. Because (ii) holds with probability at least $1/2$ and (i) with probability $1/2$, we have that the overall probability of an assignment being in the set $\sigma$ is, for some choice of target residues, at least $\frac{1}{4(k+1)B}$.

## Why it works

For any $u < v < w$, the value of $\operatorname{bit}(u, w)$ equals exactly one of $\operatorname{bit}(u, v)$ or $\operatorname{bit}(v, w)$, these two being unequal. This is because in the sequence $u \to v \to w$, the bit $\operatorname{bit}(u, w)$ must be flipped exactly once. We'll use this property several times later.

We need to show that for any configuration in $\sigma$, each prisoner knows his or her hat color with certainty. To obtain a contradiction, assume that Alice has two possible choices for her hat color: $s$ and $t$, both leading to configurations in $\sigma$: $(x_1, \ldots, x_n, y_1, \ldots, y_k)$ and $(x'_1, \ldots, x'_n, y'_1, \ldots, y'_k)$. Then, for some $z_1, \ldots, z_{k-1}$,

$x_1 = s \qquad y_1, \ldots, y_k = z_1, \ldots, z_{i-1}, \boldsymbol{t}, z_i, \ldots, z_{j-1}, z_j, \ldots, z_{k-1}$

$x'_1 = t \qquad y'_1, \ldots, y'_k = z_1, \ldots, z_{i-1}, z_i, \ldots, z_{j-1}, \boldsymbol{s}, z_j, \ldots, z_{k-1},$

If $i = j$, then the two full configurations differ by a single transposition, violating (i). So without loss of generality we can assume $i < j$. This means

$b_1, \ldots, b_{k-1} = c_1, \ldots, c_{i-2}, \boldsymbol{b_{i-1}}, \boldsymbol{b_i}, c_i, \ldots, c_{j-2}, c_{j-1}, c_j, \ldots, c_{k-2},$

$b'_1, \ldots, b'_{k-1} = c_1, \ldots, c_{i-2}, \boldsymbol{c_{i-1}}, c_i, \ldots, c_{j-2}, \boldsymbol{b'_{j-1}}, \boldsymbol{b'_j}, c_j, \ldots, c_{k-2},$

where $c_i = \operatorname{bit}(z_i, z_{i+1})$,

Because $t < z_i < z_{i+1}$, the value $b_i = \operatorname{bit}(t, z_i)$ must be different from $b'_i = \operatorname{bit}(z_i, z_{i+1})$. This means that the sequences $b_1, \ldots, b_{k-1}$ and $b'_1, \ldots, b'_{k-1}$ are different.

Despite that, the two preceding sequences are nearly identical. The only potential differences form two families: $b_{i-1}, b_i, c_{i-1}$ and $c_{j-1}, b'_{j-1}, b'_j$. For the first family, $b_{i-1} = \operatorname{bit}(z_{i-1}, t)$, $b_i = \operatorname{bit}(t, z_i)$, $c_{i-1} = \operatorname{bit}(z_{i-1}, z_i)$. Because $z_{i-1} < t < z_i$, this means that $c_{i-1}$ equals one of $b_{i-1}$ or $b_i$; let $b_*$ be the other. Similarly, $c_{j-1}$ equals one of $b'_{j-1}, b'_j$, and let $b'_*$ be the other.

This implies that $b_1, \ldots, b_{k-1}$ and $b'_1, \ldots, b'_{k-1}$ are different, yet can be made the same by removing a single element ($b_*, b'_*$) from each. But this contradicts conditions (iii) and (iv), which guarantee that the winning sequences $b_1, \ldots, b_{k-1}$ are codewords from a single-deletion-correcting code [5]. This proof was for Alice, but the same proof works for the other prisoners.

Condition (ii) reduces the range of values $b_i$ from $\lceil \log(n + k) \rceil$ to $B = \lceil 2 \log k \rceil$, which improves the density of the code from $1/O(k \log(n + k))$ to $1/O(k \log k)$. This is at the expense of disallowing sequences $b_1, \ldots, b_{k-1}$ that contain any $b_i \le b_{\max} - B$, but for $B = \lceil \log(k^2) \rceil = \lceil 2 \log k \rceil$ the birthday paradox implies that this removes only a constant fraction of sequences $b_1, \ldots, b_{k-1}$.